\documentclass[12pt]{article}
\usepackage{latexsym}
\usepackage{amsfonts}
\usepackage{amssymb}
\usepackage{amsmath}
\newtheorem{theorem}{Theorem}

\newtheorem{proposition}{Proposition}

\begin{document}
\author{Gheorghe Minea}
\title{Entropy conditions for quasilinear first order equations on nonlinear fiber bundles\\
with special emphasis on the equation of 2D flat projective structure. II.}
\maketitle
\begin{titlepage}
\begin{abstract}
We find, in an intrinsic form, a generalized Rankine-Hugoniot condition with \\respect to an entropy density that allows to give the 
proper interpretation to a \\formula of Vol'pert reducing the entropy condition on a function with bounded variation to its expression 
for generic simple jumps. It also leads to define the conservation laws only in terms of characteristics and to point out the class 
of entropy conditions coming from oriented conservation laws.

\end{abstract} 
\end{titlepage}
%\tableofcontents
\section{Introduction}
The present paper is the direct continuation of the Part I published in arXiv in December 2009 ([2]).\\
The second section is concerned with the right interpretation, in the framework of characteristics and nonlinear fibre bundles, of a 
formula of Vol'pert from [4], allowing to reduce the entropy condition on the generalized solutions with bounded variation to the 
corresponding condition on their generic simple jumps. In order of that, we found first the intrinsic concept of generalized solution 
of Rankine-Hugoniot type with respect to an \\
entropy density, as was the last one introduced in Part I of this article.\\
The third section stresses the role played by the integrability properties of the \\2-dimensional sub-bundle of sums of tangents to 
characteristics with tangents to nonlinear fibres for the question of distinguishing the entropy densities by their respective
classes of solutions. We show that a condition of complete non-integrability is necessary and sufficient for this property, while the 
complete integrability makes all classes of solutions identical and of behaviour reversible in time.\\
The fourth section deals with a special class of entropy conditions, those coming from, called by us, oriented conservation laws.
First we succeed to define and to prove local existence and unique determination of conservation laws only in terms of characteristics.
Next we identify the entropy densities defined by oriented conservation laws and prove that being a generalized Rankine-Hugoniot  
solution, with respect to such an entropy density, comes to verifying a condition of null divergence. Finally, in the case 
of a 2-dimensional space-time continuum with the completely non-integrable sub-bundle considered in the previous section, we characterize
the entropy densities coming from oriented conservation laws. Let us mention that for the equation of 2D flat projective structure
this sub-bundle is indeed completely non-integrable.

\section{Reducing the entropy condition to its\\ expression for generic simple jumps}

\S 2.1 \texttt{Intrinsic meaning of generalized Rankine-Hugoniot solutions\\ \hspace*{20pt} with respect to an entropic density}

\noindent We refer to [2] for the definitions of the entropic density $\rho$, the essentially locally bounded section $\sigma$ and the open layer 
$G$ bounding $\sigma$ on $U$. Hence $\pi: G\rightarrow U$ is a surjective submersion and $\sigma :U\rightarrow G$ is a measurable section 
of $\pi$.
\begin{equation}\label{1}
 T_{g}^{0} G:= T_{g} G_{\pi(g)} = \mathrm{ker}\,T_{g}\pi
\end{equation}
denotes the tangent to the nonlinear fiber $G_{z}$ in $g\in G_{z}$ for $z\in U$.\\
The example from [2, \S 3.5], shows that \textit{the vector bundle $T^{0} G$ over $G$ may be not\\ orientable}; in that example $U$ is 
diffeomorphic to the M\"{o}bius band and $G_{z}$ is an open connected arc in the projectivized $P(T_{z} U)$ of the tangent plane $T_{z} U$
in $z$ to $U$, continuously depending on $z\in U$.\\
We will then consider the covering with two leaves
\begin{equation}\label{2}
 p:\widetilde{U_{\pi}}\rightarrow U
\end{equation}
where $\widetilde{U_{\pi}}=\{(z,o)\arrowvert z\in U, o \,\text{orientation of} \,G_{z}\}$ (recall that $G_{z}$ are connected arcs).
Next we define the generalized function
\begin{equation}\label{3}
 RH(\rho,\sigma) : \textsl{C}_{0}^{\infty}\Gamma(\Omega(T\widetilde{U_{\pi}}))\rightarrow \textbf{R}
\end{equation}
(denoted after Rankine and Hugoniot for convenience). If $V=\mathring{V}\subseteq U$ is properly covered by $p$, so that there exists
\begin{equation}\label{4}
 o : \pi^{-1}(V)\cap G\rightarrow O(T^{0}G)
\end{equation}
continuous section, we consider also $-o,\,(-o)_{g}=-o_{g}$ and $(V,o)\subset\widetilde{U_{\pi}},\,\,
(V,-o)\subset\widetilde{U_{\pi}},\,\,\\
(V,o)\cap (V,-o)=\varnothing,\,(V,o)\cup (V,-o)= p^{-1}(V)$ (see (\ref{2})). Now for any test function 
$\zeta\in\textsl{C}_{0}^{\infty}\Gamma(\Omega(T(V,o)))$ we put
\begin{equation}\label{5}
 <RH(\rho,\sigma),\zeta> = <J(\rho,\sigma\arrowvert_{V},\tau),p_{\ast}\zeta>
\end{equation}
for $\tau$ smooth section of $\pi$ defined on $V$, with $\tau(V)\subset G$ and
\begin{equation}\label{6}
 \sigma_{z} < \tau_{z}
\end{equation}
with respect to $o$, $\forall z\in V$ (see [2, (41) - (47)] for the definition of $J(\rho,\sigma\arrowvert_{V},\tau)$). Formula (41)
from [2] shows that the definition (\ref{5}) does not depend on $\tau$ with (\ref{6}). If we consider
\begin{equation}\label{7}
 s :(V,o) \rightarrow (V,-o)
\end{equation}
canonical, we see that
\begin{equation}\label{8}
 <RH(\rho,\sigma),s^{\ast}\zeta> = -<RH(\rho,\sigma),\zeta>,
\end{equation}
$\forall \zeta\in\textsl{C}_{0}^{\infty}\Gamma(\Omega(T(V,-o)))$. Then we define $RH(\rho,\sigma)$ on $\widetilde{U_{\pi}}$ through a 
partition of unity such that (\ref{8}) holds $\forall \zeta\in\textsl{C}_{0}^{\infty}\Gamma(\Omega(T\widetilde{U_{\pi}}))$, where 
$s:\widetilde{U_{\pi}}\rightarrow\widetilde{U_{\pi}}$ has now a clear meaning.\\
In the standard case from [2, (32) - (34) and (36)], with the orientation $o$ defined by 
[2, Proposition 1], if
\begin{equation}\label{9}
 \zeta_{z} = \phi (z)\cdot\arrowvert\mathrm{d}z_{1}\wedge\dots\wedge\mathrm{d}z_{m}\arrowvert,\,\,z\in V,
\end{equation}
and $\sigma_{z} = (z,u(z))$, we get
\begin{multline}\label{10}
<RH(\rho,\sigma),\zeta> = -\int_{V}\Bigl\{\sum_{i=1}^{m} Z^{i}(z,u(z))\dfrac{\partial\varphi}{\partial z_{i}}(z) +\\
+ \Bigl[X^{m+1}(z,u(z)) + \sum_{i=1}^{m}\dfrac{\partial Z^{i}}{\partial z_{i}}(z,u(z))\Bigr]\phi(z)\Bigr\}\mathrm{d}z.
\end{multline}
We will say that $\sigma$ \textit{is a generalized Rankine-Hugoniot solution with respect to} $\rho$ \textit{on} $U$ if 
$RH(\rho,\sigma) = 0$ on $\widetilde{U_{\pi}}$, 
hence if it satisfies locally the equation in divergence form. We stress the fact that the equation has an intrinsic meaning even if
$u$, the $X^{i}$-s, the $Z^{i}$-s and $\phi$ depend on the chart. From the Corollary 1 from [2] it results that 
$I(\rho,\sigma,G)\geqslant 0\Rightarrow RH(\rho,\sigma)\geqslant 0$ and from (\ref{8}) $RH(\rho,\sigma) = 0$, that is an entropic solution 
with respect to $\rho$ is a generalized Rankine-Hugoniot solution with respect to $\rho$.\\

\noindent \S 2.2 \texttt{A formula for generalized Rankine-Hugoniot solutions\\ \hspace*{20pt} with bounded variation}

\noindent This formula is essentially due to Vol'pert [4]. Let us consider for $\sigma$ as before
\begin{equation}\label{11}
 R_{\sigma}(k)=\int_{\arrowvert k,\sigma(z)\arrowvert} T_{g}\pi\cdot\rho^{z}(\mathrm{d}g),\,\,k\in G_{z},\,z\in U;
\end{equation}
hence $R_{\sigma}(k)\in T_{z} (U),\,\forall k\in G_{z}$. Next, for $\theta\in\textsl{C}_{0}^{\infty}\Gamma(\Omega(T^{0} G))$ let
\begin{equation}\label{12}
 S(\sigma,\theta)(z) = \int_{G_{z}} R_{\sigma}(k)\cdot\theta^{z}(\mathrm{d}k),
\end{equation}
so that $S(\sigma,\theta)(z) \in T_{z} U,\,\forall z\in U$. Here 
$\rho^{z} = \rho\arrowvert_{G_{z}},\,\theta^{z} = \theta\arrowvert_{G_{z}}$.
In local coordinates, where $\rho = \lambda\otimes X$, $\lambda$ the Lebesgue measure on $\textbf{R}$, 
$\dfrac{\partial Z^{i}}{\partial y} = X^{i},\,1\leqslant i\leqslant m$, we have
\begin{equation}\label{13}
 R_{\sigma}(z,y) = \mathrm{sgn}(u(z) - y)\cdot(Z(z,u(z)) - Z(z,y)),
\end{equation}
where $Z(z,y) = (Z^{i}(z,y))_{i=1}^{m}$, and
\begin{equation}\label{14}
  S(\sigma,\theta)(z) = \int \mathrm{sgn}(u(z) - y)\cdot(Z(z,u(z)) - Z(z,y))\phi(z,y)\mathrm{d}y,
\end{equation}
if 
\begin{equation}\label{15}
\theta(z,\mathrm{d}y) = \phi(z,y)\arrowvert\mathrm{d}y\arrowvert.
\end{equation}
The Corollary from [4, \S 15.4], is the following: let $v = (v_{1},v_{2},...,v_{n}) \in BV(\Omega,\textbf{R}^{n}),\, \\
u \in BV(\Omega,\textbf{R}),\,\arrowvert v_{k}(x)\arrowvert\leqslant M\arrowvert u(x)\arrowvert,\,\forall x\in \Omega,\,
1\leqslant k\leqslant n$. Then $w:= \mathrm{sgn}(u)\cdot v$ is in $BV(\Omega,\textbf{R}^{n})$. If, moreover, $\mathrm{div} v$ is 
absolutely continuous with respect to the Lebesgue measure on $\Omega$, then for every $\phi\in \textsl{C}_{0}^{\infty}(\Omega)$
\begin{equation}\label{16}
 \int_{\Omega} \phi(x)\,\mathrm{div} w (\mathrm{d}x) = \int_{\Omega} \phi(x)\,\mathrm{sgn} (u(x))\cdot\mathrm{div} v(x)\cdot\mathrm{d}x + 
\int_{\Gamma (u)} \phi(x) <\Delta w,\nu> H_{n-1} (\mathrm{d}x),
\end{equation}
where $\Gamma (u)$ is the set of regular points of jump of $u$, $\nu$ the normal at $\Gamma (u)$ and $\Delta w = \\
=l_{\nu} w - l_{-\nu} w$ is the jump of $w$ across $\Gamma (u)$.\\
We will use this formula for our generalized function $I(\rho,\sigma,G)$ that in local coordinates reads (see [2, (37)])
\begin{multline}\label{17}
< I(\rho,\sigma,G),\psi>=\iint_{G} \mathrm{sgn}(u(z)-y)
\Bigl\{\sum_{i=1}^{m} (Z^{i}(z,u(z))-Z^{i}(z,y))\dfrac{\partial\varphi}{\partial z_{i}}(z,y)+\\
+\Bigl[X^{m+1}(z,u(z))+\sum_{i=1}^{m}\Bigl(\dfrac{\partial Z^{i}}{\partial z_{i}}(z,u(z))-\dfrac{\partial Z^{i}}{\partial z_{i}}(z,y)
\Bigr)\Bigr]\varphi(z,y)\Bigr\}\mathrm{d}z\;\mathrm{d}y, 
\end{multline}
if $\psi_{(z,y)}=\varphi(z,y)\cdot|\mathrm{d}z_{1}\wedge\dots\wedge\mathrm{d}z_{m}\wedge\mathrm{d}y|$. \\ 
For $u\in BV(U)$ generalized Rankine-Hugoniot solution with respect to $\rho$ on $U$ we apply, following Vol'pert, (\ref{16}) to 
$u(z) - y$, in 
place of $u$, and $Z(z,u(z)) - Z(z,y)$, in place of $v$, for fixed $y$. Remark that in the sense of distributions
\begin{equation}\label{18}
 \mathrm{div}_{z} (Z(z,u(z)) - Z(z,y)) = X^{m+1}(z,u(z)) + \sum_{i=1}^{m}\Bigl(\dfrac{\partial Z^{i}}{\partial z_{i}}(z,u(z))-
\dfrac{\partial Z^{i}}{\partial z_{i}}(z,y)\Bigr),
\end{equation}
hence $\mathrm{div} v$ is absolutely continuous in virtue of the hypothesis on $u$. Integrating also with respect to $y$ we get for any 
chart $\chi : V\rightarrow\textbf{R}^{n}$ on $U$
\begin{equation}\label{19}
 < I(\rho,\sigma,G),\theta\rtimes\overline{\pi^{\ast}}\chi^{\ast}(H_{m})> = 
- \int_{\chi(\Gamma(\sigma)\cap V)} <\chi_{\ast}(\Delta S(\sigma,\theta)),\nu> H_{m-1}(\mathrm{d}z),
\end{equation}
where $H_{m} = \arrowvert\mathrm{d}z_{1}\wedge\dots\wedge\mathrm{d}z_{m}\arrowvert$ (see [2, (39) and (15)]). Therefore, in 
order that a generalized Rankine-Hugoniot solution be an entropic solution it is necessary and sufficient an inequality involving the 
jump of $R_{\sigma}$ (see (\ref{11})), in each regular 
point of jump of $\sigma$, of the kind already found in the case of the simple shock in [2, (65)].

\section{Distinguishing entropy densities by their \\
respective solutions and the integrability\\
properties of the vector sub-bundle $D\oplus T^{0} F$}

We recall that (see [2, (17), (18)]) that $D$ denotes the vector sub-bundle of $TF$ of characteristic directions
\begin{equation}\label{20}
 D_{f}=\rho_{f}(T^{0}_{f}F)-\rho_{f}(T^{0}_{f}F),
\end{equation}
where the entropy density $\rho\in\textsl{C}^{1}\Gamma(\Omega(T^{0}F)\otimes TF)$ satisfies
\begin{equation}\label{21}
 T_{f}\pi\,\rho_{f}\neq 0,\;\;\forall\; f\in F.
\end{equation}
We have already considered in [2, (106)], the map
\begin{equation}\label{22}
\kappa : F\rightarrow P(TM),\; \kappa(f)=P(T_{f}\pi\;D_{f})\in P(T_{\pi(f)}M)
\end{equation}
that depends only on the characteristics and not on $\rho$.\\
In what follows a \textit{supplementary property} is imposed on \textit{the entropic solutions} (see [2, (20), (21)]):\\

(TP I) $\forall z\in V\, \exists H_{z}\subset T_{z} M,\, \mathrm{dim} T_{z} M/H_{z} = 1$, such that
\begin{equation}\label{23}
  H_{z} + T_{f}\pi\cdot D_{f} = T_{z} M,\,\,\forall f\in \lvert\sigma_{1z},\sigma_{2z}\lvert_{G_{z}}.
\end{equation}
If we denote
\begin{equation}\label{24}
 S_{f} = \rho_{f}(T_{f}^{0} F) = \mathrm{im}\,\rho_{f}
\end{equation}
the half-line with $S_{f} - S_{f} = D_{f}$, the transversality property (TP I) entails:\\

if $f_{1},\,f_{2}\in \lvert\sigma_{1z},\sigma_{2z}\lvert_{G_{z}}$ with $\kappa (f_{1}) = \kappa (f_{2})$ then $T_{f_{1}}\pi (S_{f_{1}})
= T_{f_{2}}\pi (S_{f_{2}})$.\\

As $T_{P(L)} P(V)\widetilde{=}\mathrm{Hom}(L, V/L)$ for $L$ 1-dimensional subspace of $V$, we have\\

$T_{f}(\kappa\arrowvert_{F_{z}}) : T_{f} F_{z}\rightarrow T_{\kappa(f)} P(T_{z} M)\widetilde{=}\mathrm{Hom}(T_{f}\pi\, D_{f},
T_{z} M/T_{f}\pi\, D_{f})$\\

Then there exists a subspace $\widetilde{\mathrm{im}\,T_{f}(\kappa\arrowvert_{F_{z}})}\subset T_{z} M/T_{f}\pi\, D_{f}$ such that\\

$\mathrm{im}\,T_{f}(\kappa\arrowvert_{F_{z}}) =\mathrm{Hom}(T_{f}\pi\,D_{f},\widetilde{\mathrm{im}\,T_{f}(\kappa\arrowvert_{F_{z}})})$\\

and finally there exists a subspace $\widehat{\mathrm{im}\,T_{f}(\kappa\arrowvert_{F_{z}})}\subset T_{z} M,\,
\widehat{\mathrm{im}\,T_{f}(\kappa\arrowvert_{F_{z}})}\supset T_{f}\pi\, D_{f}$ (of dimension 2 or 1)
such that
\begin{equation}\label{25}
 \mathrm{im}\,T_{f}(\kappa\arrowvert_{F_{z}})=
\mathrm{Hom}(T_{f}\pi\,D_{f},\widehat{\mathrm{im}\,T_{f}(\kappa\arrowvert_{F_{z}})}/T_{f}\pi\, D_{f}).
\end{equation}
Let $\sigma_{1}$ and $\sigma_{2}$ be two classical solutions in the neighbourhood of a point $z_{0} \in M$ with
\begin{equation}\label{26}
 \sigma_{1z_{0}}\neq \sigma_{2z_{0}}.
\end{equation}
Suppose that (TP I) holds at $z_{0}$ on $\lvert\sigma_{1z_{0}},\sigma_{2z_{0}}\lvert_{G_{z_{0}}}$ and that, moreover, it is satisfied
also the transversality condition \\

(TP II) $\exists K_{z_{0}} \subset H_{z_{0}}, \mathrm{dim}\,H_{z_{0}}/K_{z_{0}} = 1$ such that
\begin{equation}\label{27}
 K_{z_{0}} + \widehat{\mathrm{im}\,T_{f}(\kappa\arrowvert_{F_{z_{0}}})} = T_{z_{0}} M,\,\,\,
\forall f\in \lvert\sigma_{1z_{0}},\sigma_{2z_{0}}\lvert_{G_{z_{0}}}.
\end{equation}
Let $\Sigma_{0}$ be a submanifold of $M$ through $z_{0}$ of $\mathrm{codim}\,\Sigma_{0} = 2$ and $T_{z_{0}}\Sigma_{0} = K_{z_{0}}$.
It can be shown that, under these conditions, there exists a unique hipersurface $\Sigma\supset\Sigma_{0}$ (i.e. of 
$\mathrm{codim}\,\Sigma = 1$), defined in a neighbourhood $V$ of $z_{0}$ with the property that $\forall z\in\Sigma$ the subspaces
$T_{\sigma_{1z}}\pi\,D_{\sigma_{1z}}$ and $T_{\sigma_{2z}}\pi\,D_{\sigma_{2z}}$ are transversal to $T_{z}\Sigma$ in $T_{z} M$ and 
$T_{\sigma_{1z}}\pi (S_{\sigma_{1z}})$ lie on the same side of $\Sigma$, $\forall z\in\Sigma$, while $T_{\sigma_{2z}}\pi (S_{\sigma_{2z}})$ 
lie on the other side of $\Sigma,\,\forall z\in\Sigma$ (see (\ref{24})), so that the projections of the respective two characteristics 
through any $z\in V\setminus \Sigma$ differ by the fact that, as oriented by $\rho$, one of them enters in $\Sigma$ while the other 
goes out of $\Sigma$; and, finally, such that $\sigma$ defined by 
\begin{equation}\label{28}
 \sigma_{z} = \sigma_{iz},\,\,\,i\in\{1,2\},
\end{equation}
if the characteristic of $\sigma_{i}$ through $z$ enters in $\Sigma$, be an entropic solution with respect to $\rho$ in $V$.\\
We have then the following
\begin{theorem} Let $z_{0}\in M$ and $\rho$ and $m\rho$, where $m$ is a smooth and positive function, be two entropy densities on $G$
such that every local entropic solution with respect to $\rho$, defined on a neighbourhood of $z_{0}$ with values in $G$, be also a 
generalized Rankine-Hugoniot solution with respect to $m\rho$. If, moreover, in a certain $g_{0}\in G_{z_{0}}$ the function $\kappa$ 
defined by $D$ satisfies
\begin{equation}\label{29}
 T_{g_{0}}(\kappa\arrowvert_{G_{z_{0}}})\neq 0,
\end{equation}
then $m$ is constant in a neighbourhood of $g_{0}$ from $G_{z_{0}}$.
\end{theorem}
In virtue of (TP I) we can find a neighbourhood $V$ of $z_{0}$ and suitable coordinates $(x,t,y)$ on $\pi^{-1}(V)\cap G$ that bring $\rho$
to the canonical form
\begin{equation}\label{30}
 \rho = \arrowvert\mathrm{d}y\arrowvert\otimes\Bigl(\dfrac{\partial}{\partial t} + 
\sum_{i=1}^{n} X^{i} (x,t,y)\dfrac{\partial}{\partial x_{i}} +X^{n+2}(x,t,y)\dfrac{\partial}{\partial y}\Bigr),
\end{equation}
such that, if $z_{0} = (0,0)$, then
\begin{equation}\label{31}
 \kappa\arrowvert_{G_{z_{0}}}(y)=(X^{i}(0,0,y))_{i=1}^{n}\in \textbf{R}^{n}\subset P(\textbf{R}\times\textbf{R}^{n}).
\end{equation}
The consequence of Theorem 1 is that, in the hypothesis
\begin{equation}\label{32}
 T_{g}(\kappa\arrowvert_{G_{z}})\neq 0,\,\,\,\forall g\in G_{z},\,\forall z\in U,
\end{equation}
if every entropic solution with respect to $\rho$ is generalized Rankine-Hugoniot solution with respect to $m\rho$, then $m$ is constant 
on each nonlinear fiber and therefore $m\rho$ defines the same entropy condition as $\rho$ (we remarked already in [2, after (64)], that 
the entropy condition is the same when replacing $\rho$ by $m\rho$ with $m=m(z) > 0$).\\
On the other hand, it is easy to see that, in the hypothesis
\begin{equation}\label{33}
  T_{g}(\kappa\arrowvert_{G_{z}}) = 0,\,\,\,\forall g\in G_{z},\,\forall z\in U,
\end{equation}
the entropic solutions are the same for every entropy density. Indeed, in this case the local form (\ref{30}) becomes\\
\begin{equation}\label{34}
\rho = \arrowvert\mathrm{d}y\arrowvert\otimes\Bigl(\dfrac{\partial}{\partial t} + 
\sum_{i=1}^{n} X^{i} (x,t)\dfrac{\partial}{\partial x_{i}} +X^{n+2}(x,t,y)\dfrac{\partial}{\partial y}\Bigr)
\end{equation}
and the projection of the characteristics do not depend on the solution and do not meet; only the initial singularities propagate 
along them.\\
It is interesting to see how this discussion is linked to the integrability of the vector sub-bundle $D\oplus T^{0} F$ of $TF$.\\
For a general vector sub-bundle $S_{f}$ of the tangent bundle $T_{f} F$  to a manifold $F$ it is well defined a \textit{curvature tensor} 
$R_{f}\in\mathrm{Hom}(S_{f}\wedge S_{f},T_{f} F/S_{f})$ by
\begin{equation}\label{35}
 <\widetilde{\alpha_{f}},R_{f}(X_{f},Y_{f})> := <\alpha_{f},[X,Y]_{f}> = - (\mathrm{d}\alpha)_{f}(X_{f},Y_{f}),
\end{equation}
if $X$ and $Y$ are smooth vector fields and $\alpha$ is a smooth 1-form defined in a neighbourhood $W$ of $f$, such that $\forall g\in W$: 
 \begin{equation}\label{36}
 X_{g}\in S_{g},\,\,Y_{g}\in S_{g},\,\,<\alpha_{g},Z> = 0,\,\,\,\forall Z\in S_{g};
\end{equation}
in (\ref{35}) $\widetilde{\alpha_{f}}$ is naturally defined on $T_{f} F/S_{f}$. According to Frobenius theorem the sub-bundle $S$ is 
completely integrable if and only if its curvature tensor is null in every point (see, for instance, [5]).\\
We will establish now the link between the curvature tensor $R_{f}$ of $D\oplus T^{0} F$ and $T_{f}(\kappa\arrowvert_{F_{\pi(f)}})$
(see (\ref{32}) and (\ref{33})).\\
First, for an arbitrary real vector space $V$ we consider the map
\begin{equation}\label{37}
 \delta : V\setminus\{0\}\rightarrow P(V),\,\,\,\delta(v) = P(\textbf{R}\cdot v).
\end{equation}
As $T_{P(L)} P(V)\widetilde{=} \mathrm{Hom}(L,V/L)$ and $T_{v} V\widetilde{=} V$ we have 
$T_{v}\delta : V\rightarrow\mathrm{Hom}(\textbf{R}\,v,V/\textbf{R}\,v)$, more precisely
\begin{equation}\label{38}
 (T_{v}\delta\cdot w)(\lambda\,v) = \lambda\,(w + \textbf{R}\,v),\,\,\,\forall w\in V,\,\forall \lambda\in \textbf{R}.
\end{equation}
Then in virtue of the fact that $\mathrm{Ker}\,T_{v}\delta = \textbf{R}\,v$ we have the isomorphism
\begin{multline}\label{39}
\widetilde{T_{v}\delta} :V/\textbf{R}\,v\rightarrow \mathrm{Hom}(\textbf{R}\,v,V/\textbf{R}\,v),\\
\widetilde{T_{v}\delta} (w + \textbf{R}\,v)(\lambda\,v) = \lambda\,(w + \textbf{R}\,v),\,\,\forall w\in V,\,\lambda \in\textbf{R},\,
v\in V\setminus \{0\}.
\end{multline}
On the other hand we have the canonical isomorphism
\begin{equation}\label{40}
 \widetilde{T_{f}\pi} :T_{f} F/(D_{f} + T_{f}^{0} F)\rightarrow T_{\pi(f)} M/T_{f}\pi\, D_{f}.
\end{equation}
For two vector fields $X_{f}\in D_{f}\setminus \{0\},\,Y_{f}\in T_{f}^{0} F,\,f\in F$, we verify first the equality
\begin{equation}\label{41}
 T_{f}\kappa\cdot Y_{f} = T_{T_{f}\pi\cdot X_{f}} \delta\cdot T_{f}\pi\cdot [Y,X]_{f},
\end{equation}
and deduce
\begin{proposition}
If the smooth vector fields satisfy $X_{f}\in D_{f}\setminus \{0\},\,Y_{f}\in T_{f}^{0} F,\,\forall f\in F$, then for the curvature tensor
$R$ of $D\oplus T^{0} F$ we have
\begin{equation}\label{42}
 T_{f}\kappa\cdot Y_{f} = \widetilde{T_{T_{f}\pi\cdot X_{f}} \delta}\cdot\widetilde{T_{f}\pi}\cdot R_{f}(Y_{f},X_{f}).
\end{equation}
\end{proposition}
In the case of the equation of 2D flat projective structure the map $\kappa$ from (\ref{22}) is the identity of $P(TP(V))$, for $V$ real 
vector space of $\mathrm{dim} V = 3$ (see [2, \S 3.1]). Hence in that case the condition (\ref{32}) is fulfilled and then the respective
families of solutions distinguish the different entropy conditions.

\section{Entropy conditions defined by oriented \\conservation laws}

\S 4.1 \texttt{Definition of conservation laws in terms of characteristics}

\noindent A conservation law for a quasilinear equation - that we identify with its characteristics - is a vector field that locally reads
\begin{equation}\label{43}
 X = \sum_{i=1}^{m} \dfrac{\partial Z^{i}}{\partial y}(z,y)\dfrac {\partial}{\partial z_{i}} - \Bigl(\sum_{i=1}^{m}\dfrac{\partial Z^{i}}
{\partial z_{i}}(z,y)\Bigr)\,\dfrac{\partial}{\partial y}
\end{equation}
with the property that it defines the classical solutions of the equation:
\begin{equation}\label{44}
X_{(z,y)}\in D_{(z,y)},\,\,\,\forall (z,y),\,\,\,\Bigl(\dfrac{\partial Z^{i}}{\partial y}(z,y)\Bigr)_{i=1}^{m} \neq 0_{\textbf{R}^{m}},
\,\,\,\forall (z,y).
\end{equation}
It allows to write the respective quasilinear equation for classical solutions in the form
\begin{equation}\label{45}
 \sum_{i=1}^{m}\dfrac{\partial}{\partial z_{i}}[Z^{i}(z,u(z))] = 0.
\end{equation}
Let us consider the differential form
\begin{equation}\label{46}
 \alpha_{(z,y)} = 
\sum_{i=1}^{m} (-1)^{i-1}\,Z^{i}(z,y)\,\mathrm{d} z_{1}\wedge\dots\wedge\widehat{\mathrm{d} z_{i}}\wedge\dots\wedge\mathrm{d} z_{m}
\end{equation}
on the total space $F\subseteq \textbf{R}^{m}\times \textbf{R}$ of the nonlinear fibre bundle. Then	
\begin{multline}\label{47}
 \mathrm{d}\alpha_{(z,y)} = \Big(\sum_{i=1}^{m}\dfrac{\partial Z^{i}}{\partial z_{i}} (z,y)\Big)
\mathrm{d} z_{1}\wedge\dots\wedge\mathrm{d} z_{m} +\\
+ \sum_{i=1}^{m} (-1)^{i+m}\dfrac{\partial Z^{i}}{\partial y} (z,y)
\mathrm{d} z_{1}\wedge\dots\wedge\widehat{\mathrm{d} z_{i}}\wedge\dots\wedge\mathrm{d} z_{m}\wedge\mathrm{d} y
\end{multline}
or
\begin{equation}\label{48}
\mathrm{d}\alpha = i_{X}\nu,\,\,\,\,\nu = (-1)^{m-1}\mathrm{d} z_{1}\wedge\dots\wedge\mathrm{d} z_{m}\wedge\mathrm{d} y 
\end{equation}
and $X$ is given by (\ref{43}). Let in these local coordinates
\begin{equation}\label{49}
Y_{(z,y)} = \dfrac{\partial}{\partial y},
\end{equation}
a nowhere zero vector field tangent to the nonlinear fibres. From (\ref{46}) - (\ref{48}) and (\ref{44}) we infer
\begin{equation}\label{50}
 i_{Y}\alpha = 0,\,\,\,i_{X}\mathrm{d}\alpha = 0,\,\,\,\,i_{Y}\mathrm{d}\alpha_{(z,y)} \neq 0,\,\,\forall (z,y). 
\end{equation}
It is easy to see that the converse is also true, so that if $X$ is a fixed nowhere zero smooth vector field tangent to the
characteristics and $Y$ is a fixed nowhere zero smooth vector field tangent to the nonlinear fibres, a smooth $(m-1)$-form $\alpha$
with (\ref{50}) defines through a smooth nowhere zero $(m+1)$-form $\tilde{\nu}$ and the relation
\begin{equation}\label{51}
 i_{\tilde{X}}\tilde{\nu} = \mathrm{d}\alpha
\end{equation}
a smooth vector field $\tilde{X}$ nowhere zero and tangent to the characteristics (see (\ref{48})). Moreover, a smooth section
$\sigma$ of $F$ will be a classical solution of the quasilinear equation so defined by $\tilde{X}$ if and only if
\begin{equation}\label{52}
 \mathrm{d} \sigma^{\ast}\alpha =0.
\end{equation}
Indeed, a smooth section $\sigma$ is a classical solution of the quasilinear equation defined by $\tilde{X}$ if and only if
\begin{equation}\label{53}
 \sigma^{\ast}(i_{\tilde{X}}\tilde{\nu}) = 0,
\end{equation}
for a fixed nowhere zero $(m+1)$-form $\tilde{\nu}$. Therefore, from (\ref{51}), (\ref{52}) and the general equality 
$\sigma^{\ast}(\mathrm{d}\alpha) = \mathrm{d}\sigma^{\ast}\alpha$ we get (\ref{53}). The point is that \textit{the property }
$i_{Y}\alpha =0$ from (\ref{50}) \textit{ensures that} $\sigma^{\ast}\alpha$ \textit{is measurable and locally integrable even for a 
locally essentialy bounded section} $\sigma$ \textit{so that the equation} (\ref{52}) \textit{may be satisfied by such sections in the 
sense of distributions}. Moreover, the property $i_{Y}\mathrm{d}\alpha\neq 0$ ensures that the equation (\ref{52}) is nondegenerate in 
the sense that 
\begin{equation}\label{54}
 T_{f}\pi\,D_{f} \neq 0, \,\,\,\forall f
\end{equation}
(see (\ref{44}) and (\ref{22})).\\
We suppose then that (\ref{54}) holds and instead of (\ref{50}) only that
\begin{equation}\label{55}
 i_{Y}\alpha = 0,\,\,\,i_{X}\mathrm{d}\alpha = 0,\,\,\,\,(\mathrm{d}\alpha)_{g} \neq 0,\,\,\forall g\in H,
\end{equation}
if $H\subseteq F$ is the domain of $\alpha$. From (\ref{54}) we deduce that $X$ and $Y$ are linearly independent in every point and then 
(\ref{55})
implies that $(i_{Y}\mathrm{d}\alpha)_{g}\neq 0,\,\,\forall g\in H$, since $\mathrm{d}\alpha$ is a $m$-form on a $(m+1)$-dimensional 
manifold.\\ 
On the total space $(T^{\ast} M)^{\wedge p}$ of the vector bundle $\pi_{p} : (T^{\ast} M)^{\wedge p}\rightarrow M$, over a general manifold 
$M$, it is defined a canonical p-form $\theta^{p}$ as follows: first 
$$(T\pi_{p})^{\wedge p} : \{T[(T^{\ast} M)^{\wedge p}]\}^{\wedge p}\rightarrow (TM)^{\wedge p}$$
is natural; next $(T^{\ast} M)^{\wedge p}\widetilde{=} [(TM)^{\wedge p}]^{\ast}$ so that for $\alpha\in (T^{\ast} M)^{\wedge p}$ it is 
defined $\alpha\cdot (T_{\alpha}\pi_{p})^{\wedge p}\in (\{T_{\alpha}[(T^{\ast} M)^{\wedge p}]\}^{\wedge p})^{\ast}\widetilde{=}
(T_{\alpha}^{\ast}[(T^{\ast} M)^{\wedge p}])^{\wedge p}$. Then
\begin{equation}\label{56}
 \theta_{\alpha}^{p}:= \alpha\cdot (T_{\alpha}\pi_{p})^{\wedge p}.
\end{equation}
If $V = \pi(H)$, the condition $i_{Y}\alpha =0$ on $\alpha$ allows to define a bundle map
\begin{equation}\label{57}
\iota : H\rightarrow (T^{\ast}V)^{\wedge(m-1)}
\end{equation}
so that
\begin{equation}\label{58}
 \alpha = \iota^{\ast} \theta^{(m-1)},
\end{equation}
while the condition $i_{Y}\mathrm{d}\alpha\neq 0$ may be read
\begin{equation}\label{59}
 T_{f}(\iota\arrowvert_{H_{z}})\neq 0,\,\,\,\forall f\in H_{z},\,\,\forall z\in V,
\end{equation}
i.e. $\iota\arrowvert_{H_{z}}$ is an immersion in $(T^{\ast}_{z} V)^{\wedge(m-1)}$ for every $z\in V$.\\

\noindent \S 4.2 \texttt{Local existence of conservation laws for given characteristics}

\noindent Let $U$ be a manifold, $S \subset U$ a submanifold of codimension 1 in $U$ and $D$ a vector sub-bundle of $TU$ of 
$\mathrm{dim} D_{u} = 1,\,\forall u \in U$. Consider the following definition:

The triad $(U,S,D)$ will be called a \textit{neighbourhood of foliation cut} if $\exists\,p:U\rightarrow S$ surjective submersion with 
the properties:

$(i)\,\, p\arrowvert_{S} =\mathrm{id}_{S};$

$(ii)\,\,\mathrm{Ker}\,T_{u} p = D_{u},\,\,\forall u\in U;$

$(iii)\,\,p^{-1}(\{s\})$ is connected non-compact $\forall s \in S$.

Such a neighbourhood of foliation cut can always be constructed using the local flow of a smooth vector field $X$ with 
$X_{u}\in D_{u}\setminus \{0\},\,\,\forall u\in U$.
We will use the following propositions
\begin{proposition} Let $(U,S,D)$ be a neighbourhood of foliation cut and $X$ a smooth vector field on U such that 
$X_{u}\in D_{u}\setminus \{0\},\,\forall u\in U$. Let $d = \mathrm{dim}\,U$ and $p$ be such that $1\leqslant p\leqslant d-1$ be fixed. 
Then 
the system for the unknown form $\alpha\in \textsl{C}^{1}\Gamma ((T^{\ast} U)^{\wedge p})$
\begin{equation}\label{60}
(i_{X}\mathrm{d}\alpha)_{u} = \beta_{u},\, (i_{X}\alpha)_{u} =\gamma_{u},\,\forall u\in U, 
\end{equation}
with the boundary condition
\begin{equation}\label{61}
(j^{\ast}\alpha)_{s} = \delta_{s},\,\forall s\in S,
\end{equation}
where $\beta \in \textsl{C}\Gamma ((T^{\ast} U)^{\wedge p}),\,\gamma\in \textsl{C}^{1}\Gamma ((T^{\ast} U)^{\wedge(p-1)}),\,\delta \in
\textsl{C}^{1}\Gamma ((T^{\ast} S)^{\wedge p})$ are given and \\
$j:S\rightarrow U$ is the respective embedding, has a unique solution in $U$.
\end{proposition}
\begin{proposition}
 Let $D_{1}$ and $D_{2}$ be two smooth sub-bundles of $TU$ of dimension 1 and $(U_{1},S_{1},D_{1})$ and $(U_{2},S_{2},D_{2})$ be two 
neighbourhoods of foliation cut, such that $U_{2} \subseteq U_{1}\subseteq U$ and $S_{2}\subseteq S_{1}$. Then 
$\exists \,U_{1}^{\prime}\subseteq U_{2}$ such that $(U_{1}^{\prime},S_{2},D_{1})$ be a neighbourhood of foliation cut. 
\end{proposition}

The local existence and uniqueness result is contained in
\begin{theorem} Let $(U_{i},S,D_{i}),\,i=1, 2$, be two neighbourhoods of foliation cut for two \\
different foliations, with $U_{2}\subset U_{1},\;j:S\rightarrow U_{2}$ the embedding and
$X_{iu} \in D_{iu}\setminus \{0\},\\
\,\,\forall u\in U_{i},\,i=1, 2$, smooth vector fields. Then\\
$1^{0}$. $\forall\,\gamma\in\textsl{C}^{1}\Gamma ((T^{\ast} S)^{\wedge p})\,\, \exists\,
\alpha\in \textsl{C}^{1}\Gamma ((T^{\ast} U_{2})^{\wedge p})$ such that
\begin{equation}\label{62}
 (i_{X_{1}}\mathrm{d}\alpha)_{u} = 0, \,\,(i_{X_{2}}\alpha)_{u} = 0,\,\,\forall u\in U_{2},
\end{equation}
\begin{equation}\label{63}
(j^{\ast}\alpha)_{s} = \gamma_{s},\,\,\forall s\in S,
\end{equation}
of the form
\begin{equation}\label{64}
 \alpha = \beta + \mathrm{d}\phi
\end{equation}
where
\begin{equation}\label{65}
 (i_{X_{1}}\beta)_{u} = 0,\,\,\forall u\in U_{1},\,\,(j^{\ast}\beta)_{s} = \gamma_{s},\,\,\forall s\in S,
\end{equation}
and
\begin{equation}\label{66}
 (i_{X_{2}}\mathrm{d}\phi)_{u} = - (i_{X_{2}}\beta)_{u},\,\,(i_{X_{2}}\phi)_{u} = 0,\,\,\forall u\in U_{2},
(j^{\ast}\phi)_{s} = 0,\,\,\forall s\in S.
\end{equation}
Moreover, $(\mathrm{d}\alpha)_{u}\neq 0,\,\,\forall u\in U_{2}$, if and only if $(\mathrm{d}\gamma)_{s}\neq 0,\,\,\forall s\in S$.

$2^{0}$. If $\alpha\in \textsl{C}^{1}\Gamma ((T^{\ast} U_{1})^{\wedge p})$ verifies
\begin{equation}\label{67}
 (i_{X_{1}}\mathrm{d}\alpha)_{u} = 0,\,\,\forall u\in U_{1}, \,\,(i_{X_{2}}\alpha)_{u} = 0,\,\,\forall u\in U_{2},
\,\,(j^{\ast}\alpha)_{s} = 0,\,\,\forall s\in S,
\end{equation}
then $\alpha = 0$ on $U_{2}$.
\end{theorem}
The efficiency of the theorem depends on supplementary information about the two \\
foliations that would allow construction of \textit{common} neighbourhoods of foliation cut for them. An important instance is given 
by the interesting for us foliations of characteristics  and of nonlinear fibers in the case of the equation of 2D flat projective 
structure. They may also be defined as the level sets of the two submersions $\pi$ and $\mathring{\pi}$ (see [2, (71), (72) and (75)]):
\begin{equation}\label{68}
 F(V):= \{(p,d)\arrowvert\, p\in P(V),\, d\in G_{2}(V),\ p\in d\},
\end{equation}
\begin{equation}\label{69}
 \pi:F(V)\rightarrow P(V),\,\,\pi(p,d) = p;\,\,\mathring{\pi}:F(V)\rightarrow G_{2}(V),\,\,\mathring{\pi}(p,d) = d.
\end{equation} 
\begin{theorem}
For the foliations of characteristics and of nonlinear fibers of the equation of 2D flat projective structure for each 
point $f_{0} \in F(V)$ there exists a fundamental system of common neighbourhoods of foliation cut $(U_{n},S_{n},D) = 
(U_{n},S_{n},T^{0} F)$ with $f_{0}\in S_{n},\,\,\forall n$.
\end{theorem}
In fact, we consider first the domain of $\textbf{R}^{3}$,
\begin{equation}\label{70}
 U =\{(x,t,y)\arrowvert\,0 < \arrowvert t\arrowvert < x,
\,\,t^{2} + x^{2} < 1,\,\,\arrowvert y\arrowvert < 1\}
\end{equation}
and the cut $S:= \mathrm{graph} (\phi)$ where
\begin{equation}\label{71}
\phi : \{(x,t)\arrowvert\,t^{2}+ x^{2} < 1,\,0 < \arrowvert t\arrowvert < x\}\rightarrow (-1,1),\,\,\,
\phi(x,t) = -t/x.
\end{equation}  
We have $\pi(x,t,y) = (x,t),\,\,\mathring{\pi}(x,t,y) = (y,x-ty)$ and it is easy to see that the leaves of both $\pi$ and 
$\mathring{\pi}$ from $U$ are connected, non-compact and that each leaf from both families cuts $S$ in precisely one point.\\

\noindent \S 4.3 \texttt{Properties of entropy densities defined by oriented conservation laws}

\noindent Our concern here wil be the entropy densities $\rho$ for which there exists the everywhere positive density $\mu$ on $M$ (or on an 
open subset $U$ of it), and the odd form $\tau$ on $F$ (or on an open $G\subseteq F$ such that $\pi(G) = U$), that is
\begin{equation}\label{72}
 \tau = \alpha\otimes\omega,
\end{equation}
with $\alpha$ even $(m-1)$-form coming from a conservation law, i.e. with the properties (\ref{55}), and $\omega$ a local orientation on
$G$, such that 
\begin{equation}\label{73}
 \rho\rightthreetimes\overline{\pi^{\ast}}\mu = \mathrm{d}\tau.
\end{equation} 
Here $\pi$ is the usual projection $\pi:F\rightarrow M$ and for 
$\mu\in\textsl{C}^{\infty}\Gamma(\Omega(TU))$,\\
$\overline{\pi^{\ast}}\mu\in\textsl{C}^{\infty}\Gamma(\Omega(TG/T^{0}G))$ is defined by
\begin{equation}\label{74}
(\overline{\pi^{\ast}}\mu)_{g}=\mu_{\pi(g)}\circ(\overline{T_{g}\pi})^{\wedge m},
\end{equation}
where
\begin{equation}\label{75}
\overline{T_{g}\pi}:T_{g}G/T_{g}^{0}G \widetilde{\longrightarrow}T_{\pi(g)}U
\end{equation}
(see also [2,(39), (40)]). The operation appearing in the left hand side of (\ref{73}) is introduced in [2, (14), (16)], still let us 
recall it; for $W\subset V$ vector subspace there is a canonical isomorphism
\begin{equation}\label{76}
\Omega(W)\otimes\Omega(V/W)\widetilde{\longrightarrow}\Omega(V),
\end{equation}
where $\Omega(V)$ denotes the vector space of constant densities on $V$. We denoted $\rho \rtimes\sigma$ the image of
$\rho \otimes \sigma$ through this isomorphism. Next, for $\rho\in\Omega(W)\otimes V,\;\phi\in\Omega(V/W)$ we considered the contracted 
tensor product
\begin{equation}\label{77}  
\rho\rightthreetimes\phi=i_{X}(\lambda\rtimes\phi),
\end{equation} 
if $\rho=\lambda\otimes X,\;\lambda\in\Omega(W),\;X\in V$; and finally for $\rho\in\textsl{C}\Gamma(\Omega(T^{0}F)\otimes TF)$ and 
$\phi\in\textsl{C}\Gamma(\Omega(TF/T^{0}F))$ it is defined, through the contracted product (\ref{77}) taken in fibers, the odd (d-1)-form 
on $F$, if $d$ denotes $\dim F$:
\begin{equation}\label{78}
\rho\rightthreetimes\phi\in\textsl{C}\Gamma((T^{\ast}F)^{\wedge(d-1)}\otimes S(TF)).
\end{equation}
Here $S(TF)$ denotes the total space of the vector bundle of odd scalars over $F$: see [2, (1)]. \\
We call $\tau$ \textit{oriented conservation law} even if in (\ref{72}) only $\tau$ is uniquely defined; we have in mind the fact that 
a given orientation of $G$ may be only locally defined and moreover its meaning for $\tau$ is not unique: $\alpha\otimes\omega
= (-\alpha)\otimes (-\omega)$.\\
For an entropy density of the form (\ref{73}) the generalized function $I(\rho,\sigma,G)$ (see [2,(31), (37)] and (\ref{17}) here) can
be written only in terms of $\mathrm{d}\tau$. In order of that, 
for $\zeta\in\textsl{C}_{0}^{\infty}\Gamma(\Omega(T^{0} G))$ we consider its fiber primitive 
$\int\zeta\in\textsl{C}^{\infty}(\tilde{G},\textbf{R})$ defined by
\begin{equation}\label{79}
 \Bigl(\int\zeta\Bigr)(g,o) := \int_{(-\infty,\;g)_{o}} \zeta\in \textbf{R}.
\end{equation}
Recall that
\begin{equation}\label{80}
 \varpi:\tilde{G}\rightarrow G,\,\,\,\tilde{G}_{g}=\varpi^{-1}(\{g\}):=O(T_{g}G_{\pi(g)})
\end{equation}
is a natural covering of $G$ with two leaves ($O(V)$ denotes the set of two orientations of the vector space $V$) and
\begin{equation}\label{81}
 (-\infty,\;g)_{o}=\{h\in G_{\pi(g)}\arrowvert \;h<g \;with\; respect\; to\; o\}.
\end{equation}
In [2, (27)] we gave an analogous definition for the fiber primitive
\begin{equation}\label{82}
 \int\psi\in\textsl{C}^{\infty}\Gamma(\Omega(T\tilde{G}/T^{0}\tilde{G}))
\end{equation}
when $\psi\in\textsl{C}_{0}^{\infty}\Gamma(\Omega(TG))$. If we replace in the definition  
\begin{equation}\label{83}
< I(\rho,\sigma,G),\psi>= \int_{ D(\sigma)}\mathrm{d}(\varpi^{\ast}\rho\rightthreetimes\int\psi).
\end{equation}
(see [2, (31)]) $\psi = \zeta\rtimes\overline{\pi^{\ast}}\mu$, from (\ref{73}) we get
$$ \varpi^{\ast}\rho\rightthreetimes\int\psi = \Bigl(\int\zeta\Bigr)\cdot\varpi^{\ast}(\mathrm{d}\tau).$$
Thus, in the case of an entropy density coming from an oriented conservation law through (\ref{73}), it is natural to consider the 
generalized function
\begin{equation}\label{84}
 <I(\tau,\sigma,G),\zeta> :=\int_{ D(\sigma)}\mathrm{d}\Bigl[\Bigl(\int\zeta\Bigr)\cdot\mathrm{d}(\varpi^{\ast}\tau)\Bigr],
\end{equation}
for $\zeta\in\textsl{C}_{0}^{\infty}\Gamma(\Omega(T^{0} G))$. Analogously, in place of $< J(\rho,\sigma\arrowvert_{V},\tau),\zeta>$ 
(see [2, (41)]), we will consider the distribution
\begin{equation}\label{85}
 < J(\tau,\sigma\arrowvert_{V},\eta),\phi> :=\int_{\lvert\sigma,\eta\lvert}\mathrm{d}[(\phi\circ\pi)\cdot\mathrm{d}\tau]
+\int_{(\mathrm{im}\;\eta,\;o\;(\eta\rightarrow\sigma))}(\phi\circ\pi)\cdot\mathrm{d}\tau,
\end{equation}
for $\phi\in\textsl{C}_{0}^{\infty}(V,\textbf{R})$. Accordingly, we will define the distribution (compare with (\ref{3}),(\ref{5}))
\begin{equation}\label{86}
 RH(\tau,\sigma) : \textsl{C}_{0}^{\infty}(\widetilde{U_{\pi}},\textbf{R})\rightarrow \textbf{R}
\end{equation}
by
\begin{equation}\label{87}
 <RH(\tau,\sigma),\phi> = <J(\tau,\sigma\arrowvert_{V},\eta),\phi\circ p^{-1}>
\end{equation}
for $\phi\in\textsl{C}_{0}^{\infty}((V,o),\textbf{R})$; we keep in mind that $p: (V,o)\rightarrow V$ is a diffeomorphism (see
the notation $(V,o)$ after (\ref{4})). \\
For $W\subset V$ vector subspace, $\mathrm{dim}\,V = d,\,\,\mathrm{dim}\,W = p$ and $\pi:V\rightarrow V/W$ canonical, there is a natural 
isomorphism
\begin{equation}\label{88}
O(W)\otimes O(V/W)\widetilde{\longrightarrow} O(V),\,\, o\otimes \vartheta\mapsto o\rtimes\vartheta,
\end{equation}
\begin{equation}\label{89} 
(o \rtimes\vartheta)(v_{1}\wedge\dots\wedge v_{d})= o(Pv_{1}\wedge\dots\wedge Pv_{p})\vartheta(\pi v_{p+1} \wedge\dots\wedge\pi v_{d})        
\end{equation}
if $\pi v_{p+1} \wedge\dots\wedge\pi v_{d}\neq 0$ and $P: V\rightarrow W$ is the projection corresponding to the direct sum decomposition 
$$V=W\dotplus\sum_{j=p+1}^{d}\textbf{R}\cdot v_{j}.$$ 
Analogously, there is the isomorphism
\begin{equation}\label{90}
O(V/W)\otimes O(W)\widetilde{\longrightarrow} O(V),\,\, \vartheta \otimes o\mapsto \vartheta\ltimes o,
\end{equation}
\begin{equation}\label{91}
(\vartheta\ltimes o)(v_{1}\wedge\dots\wedge v_{d})=\vartheta(\pi v_{1} \wedge\dots\wedge\pi v_{d-p})o(Pv_{d-p+1}\wedge\dots\wedge Pv_{d}) 
\end{equation}
if $\pi v_{1} \wedge\dots\wedge\pi v_{d-p}\neq 0$ and $P: V\rightarrow W$ is the projection corresponding to the direct sum decomposition 
$$V=W\dotplus\sum_{j=1}^{d-p}\textbf{R}\cdot v_{j}.$$ 
We have always
\begin{equation}\label{92}
 \vartheta\ltimes o = (-1)^{p(d-p)}o \rtimes\vartheta.
\end{equation}
If $o\in O(W)$ and $\omega\in O(V)$ we denote $\omega\diagup o\in O(V/w)$, respectively $o\diagdown\omega\in O(V/W)$, the
orientations verifying
\begin{equation}\label{93}
 o\rtimes (\omega\diagup o) =\omega,\,\,\,(o\diagdown\omega)\ltimes o =\omega.
\end{equation}
Of course
\begin{equation}\label{94}
 o\diagdown\omega =(-1)^{p(d-p)} \omega\diagup o.
\end{equation}
If $G\rightarrow U$ is an open layer bounding the locally essentialy bounded section $\sigma$, defined on $U$, in the sense that 
(see [2, (21)]): $\forall\;z_{0}\in U\; \exists V=\mathring{V}\subseteq U,\;V\ni z_{0},
\;\exists \;\sigma_{1},\;\sigma_{2}\in\textsl{C}\Gamma(G\arrowvert V)$ such that
\begin{equation}\label{95}
\sigma_{z}\in \lvert\sigma_{1z},\sigma_{2z}\lvert_{G_{z}},\;z\in V,
\end{equation}
and $\omega$ is a local orientation on $G$ and $o$ a local orientation of the nonlinear fibers of $G$, we will consider the local 
orientations $o\diagdown\omega$ and $\omega\diagup o$ of the vector bundle $TG/T^{0} G$ defined in the fibers as before.\\
The definition of the entropy density coming from an oriented conservation law (\ref{73})
is justified by the following
\begin{theorem} Let $\tau$ be an odd $(m-1)$-form on $G\subseteq F$, such that the form $\alpha$ defined locally by (\ref{72}) satisfies 
(\ref{55}). Let $\iota$ be defined by (\ref{57}) and (\ref{58}) if $H\subseteq G$ is orientable by $\omega$, so that $\alpha$ is itself 
defined on $H$.\\
Let $\sigma$ be a locally essentially bounded section of $F$, defined on $U = \pi(G)$, whose image is bounded by the open layer $G$.\\  
Then $\sigma^{\ast}\alpha$ is correctly defined on $V:= \pi(H)$ by
\begin{equation}\label{96}
\sigma^{\ast}\alpha = \iota\circ\sigma
\end{equation}
as this relation holds for $\sigma$ differentiable. Moreover $\sigma\circ p: \widetilde{U_{\pi}}\rightarrow G$ is canonically oriented
by the correspondence of orientations
\begin{equation}\label{97}
 ((\sigma\circ p)^{\ast}\omega)_{(z,o)} = [(T_{(z,o)} p)^{-1}\cdot\overline{T_{\sigma(z)}\pi}]_{\ast}(\omega\diagup o)_{\sigma(z)}
\end{equation}
for $\omega$ local orientation on $G$, since the orientation 
$[(T_{(z,o)} p)^{-1}\cdot\overline{T_{g}\pi}]_{\ast}(\omega\diagup o)_{g}\in O(T_{(z,o)}\widetilde{U_{\pi}})$
is well defined and independent of $g$, with $\pi(g) = z$. Therefore the pull-back
\begin{equation}\label{98}
(\sigma\circ p)^{\ast}\tau = p^{\ast}(\sigma^{\ast}\alpha)\otimes (\sigma\circ p)^{\ast}\omega
\end{equation}
is well defined, locally essentially bounded and independent of the local representation (\ref{72}) of $\tau$. And finally
\begin{equation}\label{99}
 RH(\tau,\sigma) = -\mathrm{d}((\sigma\circ p)^{\ast}\tau)
\end{equation}
where the exterior differentiation is taken in the sense of distributions.
\end{theorem}
In the computations above we used the calculus with odd differential forms from [3]. Remark that in the right hand side of (\ref{99}) there
is a generalized odd differential $m$-form on the $m$-dimensional manifold $\widetilde{U_{\pi}}$, hence a generalized density, or 
distribution, like in the left hand side.\\
Taking into account the importance of the expressions $R_{\sigma}$ and $S(\sigma,\theta)$ (see (\ref{11}) and (\ref{12})) in the formula
(\ref{19}), it is of interest their form for the special case of an entropy density defined by an oriented conservation law.
\begin{theorem} Let the entropy density $\rho$ be defined by (\ref{73}) and (\ref{72}) on an open, oriented by $\omega$, subset $H$ of $F$
and $\iota$ the immersion defined, in (\ref{57}) and (\ref{58}), by $\alpha$ submitted to (\ref{55}). If, in addition, $T^{0}H$ is oriented 
by $o$ we have
\begin{equation}\label{100}
  i_{\int_{\arrowvert a,b\arrowvert} T_{g}\pi\cdot\rho^{z}(\mathrm{d} g)} \mu_{z} = 
\mathrm{sgn} (b - a)(\iota (b) - \iota (a)) \otimes (- \omega \diagup o),
\end{equation}
where
\begin{center}
$\mathrm{sgn} (b - a) =\Biggl\{
\begin{array}{lll}
 1,\; b > a,\\
 0,\; b = a,\\
-1,\;b < a,
\end{array}$
with \;respect \;to $o$ on $H_{z}$,
\end{center}
the difference $\iota (b) - \iota (a)$ is taken in $(T^{\ast}_{z} M)^{\wedge (m-1)}$, the orientation $- \omega \diagup o$ is thought
on $V:=\pi(H)$ and $z\in V$. In other words (see (\ref{89}) and (\ref{93})):
\begin{multline}\label{101}
\mu_{z}(\int_{\arrowvert a,b\arrowvert} T_{g}\pi\cdot\rho^{z}(\mathrm{d} g),v_{1},v_{2},\dots,v_{m-1}) = \\
= - \mathrm{sgn} (b - a)(\iota (b)_{z} -\iota (a)_{z})(v_{1},v_{2},\dots,v_{m-1})
\omega_{h}(Y\wedge W\wedge V_{1}\wedge V_{2}\wedge\dots\wedge V_{m-1})/o_{h}(Y),
\end{multline}
if $\pi(h) = z,\,\,Y \in T_{h} H_{z},\,\,W\in T_{h} H,\,\, V_{i} \in T_{h} H,\,\,v_{i} = T_{h}\pi\,\,V_{i},\,i = 1,\dots, m-1,\\
\int_{\arrowvert a,b\arrowvert} T_{g}\pi\cdot\rho^{z}(\mathrm{d} g) = T_{h}\pi\,\,W,\,\,
Y\wedge W \wedge V_{1}\wedge V_{2}\wedge\dots\wedge V_{m-1} \neq 0$.
\end{theorem}
Remember that the definition of an entropy density $\rho$ defined by an oriented conservation law $\tau$ is linked to the 
existence of an everywhere positive density $\mu$ on the base space-time continuum, such that the equality (\ref{73}) holds. 
This property is clarly invariant at multiplication by an everywhere positive function on the base manifold: $f\rho$ verifies 
the same equality using $(1/f)\;\mu$ in place of the needed density, when $f = f(z) >0$. That verification was necessary anyhow, 
because we knew that $f \rho$ defines the same entropy condition as $\rho$. Then the following result has its own interest.
\begin{proposition} An everywhere positive smooth function $f$ defined on the total manifold $F$, or only on an open subset of it, 
has the property that for every entropy density $\rho$ of the form (\ref{73}) $f\rho$ is of the same form if and only if $\exists\; g > 0,
\;h > 0$ such that
\begin{equation}\label{102}
 f = g\cdot h,\;\mathcal{L}_{X} g = 0,\;\mathcal{L}_{Y} h = 0,
\end{equation}
for vector fields $X,\;Y$ such that $X_{e}\in D_{e}\setminus\{0\},\;Y_{e}\in T^{0}_{e}F\setminus\{0\},\;\forall e$.
\end{proposition}
Here $\mathcal{L}_{X}\phi = <\mathrm{d}\phi,X>$ denotes the usual derivative of the function $\phi$ along $X$.\\

\noindent \S 4.4 \texttt{Characterization of entropy densities defined by oriented conservation\\ \hspace*{18pt} laws for a base manifold of dimension 2 and completely non-integrable\\ \hspace*{18pt} sub-bundle $D\oplus T^{0} F$}.

\noindent Suppose that $\mathrm{dim} M = 2$ and $R_{e} \neq 0,\;\forall e$ in the common domain of definition of the smooth vector fields $X$ 
and $Y$, with $X_{e}\in D_{e}\setminus \{0\},\,Y_{e}\in T_{e}^{0} F,\,\forall e$ (see Section 3 here, especially (\ref{42})). As we 
stressed already, the equation of 2D flat projective structure satisfies this condition. These two conditions are equivalent to
\begin{equation}\label{103}
\textbf{R} X_{e} + \textbf{R} Y_{e} + \textbf{R} [X,Y]_{e} = T_{e} F,\;\forall e.
\end{equation}
We have the following nice complement of the Proposition 3 above
\begin{proposition}
Under the condition (\ref{103}) a positive differentiable function $f$, defined on a common simply connected domain with $X$ and $Y$,
is of the form (\ref{102}) if and only if the differential 1-form $\alpha$ defined by
\begin{equation}\label{104}
 <\alpha,X> = \mathcal{L}_{X} \mathrm{log} f,\;<\alpha,Y> = 0,\;<\alpha,[Y,X]> = \mathcal{L}_{Y}\mathcal{L}_{X}\mathrm{log} f,
\end{equation}
is closed: $\mathrm{d}\alpha = 0$.
\end{proposition}
Let us recall the following definition: the divergence of a vector field $X$ with respect to a smooth nowhere null density $\mu$ on the 
same manifold is the function that verifies $\mathrm{div}_{\mu} (X)\cdot\mu = \mathcal{L}_{X}\mu$, where the right hand side represents 
the Lie derivative of $\mu$ with respect to $X$.  We have also $\mathcal{L}_{X}\mu =\mathrm{d} (i_{X}\mu)$.
Then the following characterization holds:
\begin{theorem}
 Let the condition (\ref{103}) be satisfied and $\rho$ be an entropy density of the form $\rho =\lambda\otimes X$. Then there exists 
$f$ smooth and everywhere positive on a simply connected subset $V$ of the projection of the domain of $\rho$ on the base manifold $M$, 
such that   
\begin{equation}\label{105}
\mathrm{d}(f\;\rho\rightthreetimes\overline{\pi^{\ast}}\mu) = 0 
\end{equation}
if and only if, for $\nu := \lambda\rtimes\overline{\pi^{\ast}}\mu$, the 1-form $\xi$ defined by
\begin{equation}\label{106}
 <\xi,X> = \mathrm{div}_{\nu}(X),\;<\xi,Y>= 0,\; <\xi,[Y,X]> = \mathcal{L}_{Y}\mathrm{div}_{\nu}(X),
\end{equation}
on $\pi^{-1} (V)$, is closed: $\mathrm{d}\xi = 0$.
\end{theorem}

\vspace*{20pt}
\noindent {\large {\bf Aknowledgement.}}  Support from the Grant PN2 No. 573/2009.

\end{document}